\documentclass[11pt]{article}

\usepackage{amsfonts}
\usepackage{amssymb}
\usepackage{amsmath}
\usepackage{amsthm}

\newtheorem{definition}{Definition}

\newtheorem{theorem}{Theorem}

\newtheorem{lemma}{\noindent Lemma}

\newtheorem{problem}{\noindent Problem}

\def\Ab{\operatorname{Ab}}
\def\Mar{\operatorname{Mar}}
\def\Prob{\operatorname{Prob}}
\def\Meas{\operatorname{Meas}}

\title{A method of defining central and Gibbs measures and the ergodic method}
\author{A.~M.~Vershik\footnote{St.~Petersburg Department of Steklov Institute of Mathematics, St.~Petersburg State University, and Institute for Information Transmission Problems.}}
\date{}

\begin{document}
\maketitle


\begin{abstract}
We formulate a general statement of the problem of defining invariant measures with certain properties and suggest an ergodic method of perturbations for describing such measures.
\end{abstract}


\section{Defining probability measures via projective conditional measures}

In the second half of the last century, a new method of defining probability measures in infinite-dimensional systems, alternative to the classical (Kolmogorov's) one, was gradually developed. Instead of a system of consistent finite-dimensional distributions, which  uniquely defines a measure via projections, the new method involves another data system, roughly speaking, a consistent system of conditional measures. The method appeared independently in the theory of Markov processes (E.~B.~Dynkin), in statistical physics (R.~L.~Dobrushin), etc. We present an abstract version of the method, regarding it simultaneously as a far-reaching generalization of the theory of measurable partitions of Lebesgue spaces and Rokhlin's systems of conditional measures and as the problem of describing invariant measures in the theory of dynamical systems and Gibbs measures. The presentation is in terms of {\it equipped equivalence relations} ({\it e.e.r.}), or, in other words, Borel partitions in a standard Borel space and ``projective conditional measures'' on elements of these partitions. We could also use the language of groupoids or the language of the theory of extensions of measures from special algebras (but not $\sigma$-algebras) of sets to a
 $\sigma$-algebra (not the whole $\sigma$-algebra) of measurable sets. Similar considerations in smaller generality can be found in~\cite{Ke,Sch}.

The new method has essentially introduced a large number of combinatorial, analytic, and algebraic problems: about central measures on path spaces of graded graphs, about Markov measures with given cotransitions, etc.

Consider a standard Borel space $(X, \mathfrak A)$ (isomorphic to the interval $[0,1]$ regarded as a Borel space) with the $\sigma$-algebra~${\mathfrak A}$ of all Borel sets. Let $\tau$~be a Borel equivalence relation (\textit{e.r.}) (i.e., partition) on
this space with countable classes and $\rho$~be a Borel $2$-cocycle with nonnegative real values on this relation, i.e., a Borel function on pairs of equivalent points satisfying the conditions
 $\rho(x,x)=\rho(x,y)\rho(y,x)=1$, $\rho(x,y)\rho(y,z)=\rho(x,z)$. The cocycle~$\rho$ uniquely (up to a positive coefficient) defines a finite or $\sigma$-finite nonnegative measure (in short, ``conditional projective measure'') on each equivalence class.
The pair $(\tau,\rho)$ will be called an {\it equipped equivalence relation} on the space $(X, \mathfrak A)$.\footnote{It would be natural, extending Rokhlin's terminology, to introduce the term ``semimeasurable partition with a system of conditional projective measures.''} If all classes are finite, then the e.e.r.\ is nothing else than a Borel measurable partition, and the cocycle defines conditional probability measures on all classes.

On the other hand, it is well known that if $\mu$~is a probability measure on~$X$ (i.e., on a Lebesgue space) and  $\tau$ is an e.r., then $\mu$~uniquely defines an equipment of this e.r., that is, a $2$-cocycle, or a conditional projective measure on almost every equivalence class. In the case where the e.r.\ is the partition into the orbits of an action of a countable group with a quasi-invariant measure, this cocycle is the so-called Radon--Nikodym cocycle~$RN_{\mu}$. If the cocycle is identically equal to~$1$, then the measure is said to be invariant. If the cocycle is not indicated, it is assumed to be identically equal to~$1$. We now formulate the main problem.

\begin{problem} Let $(\tau,\rho)$ be an e.e.r.\ on a standard Borel space $(X,\mathfrak A)$; find all Borel probability measures~$\mu$ for which the Radon--Nikodym cocycle~$RN_{\mu}$ coincides with~$\rho$ almost everywhere with respect to the measure~$\mu$, or, in other words, find all probability measures with given conditional projective measures on equivalence classes.
\end{problem}

In the case where such a measure is unique (it is this case that corresponds to Kolmogorov's system of finite-dimensional distributions), we may speak of extending the measure to the $\sigma$-algebra of all measurable sets; in the general case, there can be no uniqueness.

The family of all measures on $(X,\mathfrak A)$ determined by Problem~1 is well defined and forms, in a natural way, a Choquet simplex. The set of its extreme points (Choquet boundary) is called the
{\it absolute of the e.e.r.~$(\tau,\rho)$} and denoted by~$\Ab(X,\tau, \rho)$.
Any two different measures from $\Ab(X,\tau, \rho)$ are mutually singular and defined on different complete $\sigma$-algebras.

The traditional construction of Gibbs measures, as well as the problem of describing invariant measures of group actions in dynamics, obviously, fits in this scheme. The notion of absolute is closely related to various notions of boundary.

If all classes of an e.e.r.\ are finite, then we have a well-defined Borel quotient space~$X/{\tau}$, which coincides with $\Ab(X,\tau, \rho)$, and describing all (not ergodic) measures from
${\cal M}_{(\tau,\rho)}$ reduces to indicating a measure on this quotient space. But if the e.r.\ does not define a Borel quotient space (space of classes), then studying the structure of the absolute is a difficult problem and depends significantly on the geometry of classes. The solution of Problem~1 can be ``wild,'' i.e., it may happen that the absolute has no reasonable parametrization, but in many  (for example, combinatorial) problems, a parametrization can be found.

It is not difficult to extend all these definitions to e.r.'s~$\tau$ whose equivalence classes are not countable, but endowed with a well-defined locally compact topology.

It is important to emphasize that an e.r.\ can be studied only together with a cocycle, i.e., a system of conditional projective measures (even if the cocycle is identically equal to~$1$). The main role in the further analysis of the subject (uniqueness, special properties of measures, etc.)\ must be played by the geometry of classes of e.e.r.'s, but it is still poorly studied.

We now state the inverse problem.

\begin{problem} Let $M$ be a family of probability measures defined on the $\sigma$\nobreakdash-algebra~$\mathfrak A$  of a standard Borel space~$X$.  Find the minimal e.r.~$\tau$ for which all measures $\mu\in M$ define the same cocycle
$\rho\equiv RN_{\mu}$.
\end{problem}

This problem is a generalization of the traditional problem related to sufficient statistics (cf.~\cite{DF}), in which one usually seeks only measurable (for example, finite) equivalence relations. In the above setting, there are no restrictions on the e.r. For example, consider the set of all Bernoulli measures
$\prod_1^{\infty}(p,1-p)$, $p\in [0,1]$, on the space of sequences
$\prod_1^{\infty}\{0;1\}$. The required e.e.r.\ is the de Finetti partition with cocycle identically equal to~$1$: two sequences are equivalent if they coincide from some index~$n$ on and have the same number of zeros among the first $n$~coordinates.

\section{Hyperfinite and tame e.e.r.'s, the universal Markov model}

We now consider Problem~1 for a special case whose importance is due to a~large number of applications. Namely, it includes the problem of describing characters of locally finite groups, or, more generally, describing traces on AF-algebras, as well as the problem of describing central measures on path spaces of graded graphs.

{\it An equipped equivalence relation~$\tau$  is said to be hyperfinite if it is a~monotonely increasing limit of a sequence of finite equivalence relations:
$\tau=\bigcup_n  \xi_n$.} Thus, a hyperfinite e.e.r.\ can be defined by a sequence of its finite approximations, i.e., a decreasing sequence of measurable partitions $\{\xi_n\}_n$ with finite elements and conditional measures on them. Such sequences are called {\it filtrations}. For details, see~\cite{V17}.

By a number of well-known theorems, the orbit partition for a group action with an invariant measure is hyperfinite if and only if the group is amenable. However, orbit partitions with nontrivial cocycles can be hyperfinite also for nonamenable groups. Note that the hyperfiniteness condition for an e.e.r.\ is a condition on the cocycle, i.e., on the conditional projective measures, but apparently it has not been stated in this form. For Lebesgue spaces, a hyperfinite e.e.r.\ is unique up to isomorphism (generalized H.~Dye's theorem).

We impose a slightly stronger (than hyperfiniteness) condition on the approximating sequence of measurable partitions $\{\xi_n\}$:  an e.r.\ is said to be {\it tame}, or {\it locally hyperfinite}, if for every~$n$ the number of types of conditional measures of the partition~$\xi_n$  is finite. This condition singles out a class of hyperfinite e.e.r.'s that is of most interest for applications. Now we describe a universal model of a tame e.e.r.

\begin{definition}
Let $X_n$ be a finite or compact space and $\{\pi_n\}$ be a set of ``transition operators'' that send a point
$x\in X_n$ to a subset $\pi_n(x)\subset X_{n+1}$. The corresponding Markov (nonstationary) compactum~$\Mar$ is the space of sequences
$$
\Mar \subset \{\{x_n\}_{n=1}^{\infty}: x_n\in X_n,\; n=1,2, \dots\}
$$
where
 $$
 \{x_n\}\in \Mar \iff  x_{n+1}\in \pi_n(x_n)\quad\text{for every $n\geq 1$}.
 $$
\end{definition}

Elements of $\Mar$ are called trajectories, or paths. The {\it tail equivalence relation}~$\tau$ on the Markov compactum~$\Mar$ is the following relation on trajectories:
$$\{x_n\}\sim_\tau \{y_n\} \iff \quad \text{there exists $N$ such that $x_n=y_n$ for every $n>N$}.
$$

The Markov compactum $\Mar$ is endowed with the weak topology and Borel structure. A Markov Borel measure~$P$ is defined by an initial distribution~$\mu_1(\cdot)$
of the coordinate~$x_1$ and a collection of transition probabilities, i.e., a family of measures
 $\{P_{n,x}\}$, $n=1,2,\dots$, $x\in X_{n}$, where $P_{n,x}(y)=\Prob(x_{n+1}=y|x_n=x)$.

But we will need another  data system on a Markov compactum, a~{\it system of cotransition probabilities}. It is a family of measures~$\{P^{n,x}\}$ on $X_n$, ${n=1,2,\dots}$, $x\in X_{n+1}$,
where $P^{n,x}(y)=\Prob(x_n=y|x_{n+1}=x)$. Such a system does not yet define a global measure on the entire Markov compactum.

\begin{lemma}
Every system of cotransition measures defines a cocycle on the tail e.r.\ of the Markov compactum: the quotient of the conditional measures of two paths $\{x_n\}$ and $\{y_n\}$ coinciding for $n>N$ is equal to the quotient of the products of the corresponding transition probabilities
$$\prod_{1\leq i\leq N} \frac{\Prob(x_i|x_{i+1})}{\Prob(y_i|y_{i+1})}.
$$
\end{lemma}

Such cocycles will be called Markov cocycles; a Markov compactum equipped with a Markov cocycle, i.e., a system of cotransitions, will be called an {\it equipped Markov compactum}.

Obviously, every Markov measure on~$\Mar$ uniquely determines a Markov cocycle, but, in general, a system of contransitions, i.e., a Markov cocycle, does not uniquely determine a Markov measure. It is equally clear that on a Markov compactum there can exist non-Markov cocycles.

\begin{theorem}[universal model]
For every standard Borel space~$X$ and a~tame equipped equivalence relation~$\tau$ on~$X$
with cocycle~$\rho$ there exists an equipped Markov compactum~$\Mar$ and a Borel isomorphism $T:X\rightarrow \Mar$ that sends $\tau$ to the tail e.r.\ on~$\Mar$ and sends~$\rho$ to a Markov cocycle.

Thus, Problem~{\rm1} of finding all invariant measures for a tame e.e.r.\ reduces to the problem of finding all Markov probability measures~$P$ on a compactum~$\Mar$ with a given system of cotransition probabilities. In other words, to the problem of describing Markov chains with given cotransitions.
\end{theorem}

If the cocycle is identically equal to~$1$, i.e., all conditional measures of all orders are uniform, then we obtain the problem of describing all measures of maximal entropy on a given Markov compactum.

The absolute of a Markov compactum~$\Mar$ is denoted by~$\Ab(\Mar)$. The proof of the theorem essentially follows from the results of~\cite{V17}.

Instead of the language of Markov compacta, one can use the language of $\mathbb N$-graded graphs (Bratteli diagrams): the path space of such a graph is a Markov compactum, which ensures the parallelism.  In many situations (mainly of combinatorial nature), the language of graphs is more preferable. The author does not know models similar to the universal model for general e.e.r.'s.

\section{The ergodic method and the method of perturbations}

By the ergodic method of solving Problem~1 about invariant measures for hyperfinite equivalence relations we mean the method of finding invariant distributions and invariant measures based on the pointwise ergodic theorem or, more exactly, on the pointwise martingale convergence theorem, applied to the characteristic functions of sets from a basis of the $\sigma$-algebra on which the required measure is defined. In this meaning, the term is used in the author's paper~\cite{V1} and in earlier papers (see, e.g.,~\cite{Lin}). But the practical task of finding invariant measures, i.e., probabilities of cylinders, or transition probabilities, as limits of some conditional expectations can be quite difficult. A very important factor is the choice of a basis of sets whose measures are being calculated. But, on the other hand, the problem of finding all ergodic measures can be ``wild,'' so complete calculations can be essentially infeasible. A reasonable classification of hyperfinite absolutes (i.e., systems of conditional projective measures) is hardly possible; while the Borel classification of e.e.r.'s is, on the contrary, too rough (see~\cite{Ke}); the author does not know any other, intermediate, classification criteria. That is why, it is important to have feasible solvability criteria for the problem of describing invariant measures, as well as methods of reduction of problems to a few canonical problems.

One of these fundamental problems, whose solution is obtained by a~canonical application of the ergodic method, is the problem of describing all ergodic measures on the infinite product
$X^{\infty}=\prod_1^{\infty}X$ (where $X$~is a Borel space) invariant under the group~$S_{\infty}$ of all finite permutations of coordinates. Denote by~$\tau^F$ the e.r.\ on~$X^{\infty}$ generated by the partition into the orbits of the action of~$S_{\infty}$. The answer to Problem~1 is given by de Finetti's theorem and says that every ergodic invariant measure is a Bernoulli measure with an arbitrary one-dimensional distribution (= measure on~$X$). Thus, $\Ab(X^{\infty},\tau^F)=\Meas(X)$.

If we regard this answer up to a measure-theoretic isomorphism, it turns out that the absolute consists of a unique purely continuous measure on~$X$, a continuum of discrete measures, and their mixtures, i.e.,
$\Ab(X^{\infty},\tau^F) =\{\{\alpha_n\}:\alpha_n\ge0,\,\sum_n\alpha_n\leq 1\}$.

{\it In many recent examples of problems related to absolute in combinatorial and algebraic situations, the answer (hypothetical or otherwise) has a~similar structure: the absolute is a Choquet simplex (which can be called the secondary simplex), i.e., the ergodic measures themselves also admit a decomposition.} That is why, it is natural to conjecture that a proof of this fact should be sought not in a direct calculation but in analyzing the reduction to the fundamental de Finetti problem described above. We suggest the following method, which can be called the method of perturbations. For an unperturbed problem, we take the de Finetti problem about~$\tau^F$. The first step is to construct a homomorphism~$T$ of the space where we seek the absolute for some e.r.~$\tau$ to the space~$X^{\infty}$ such that  $T(\tau)$ is a subpartition of~$\tau^F$. At the second step, we must verify that the absolute of~$T(\tau)$ is similar to, or even coincides with, the absolute of~$\tau^F$.
To find~$T$, if it exists, is the most nontrivial part of the method.

The second step is related to a problem about the infinite product~$X^\infty$ which is of independent interest.

\begin{problem}
For what e.r.'s $\tau$ satisfying the condition $\tau \succ \tau^F$  does the absolute $\Ab(X^{\infty}, \tau)$
consist of all Bernoulli measures?
\end{problem}

The following partial answer is useful:

\begin{lemma} Let $\tau, \tau'$ be two e.e.r.'s on a Borel space~$X$ with cocycles identically equal to~$1$, and let
 $\Ab(X,\tau')\subset \Ab(X,\tau)$. These absolutes coincide if and only if the following condition is satisfied: for every ergodic measure  ${\mu\in \Ab(X,\tau')}$, the e.e.r.~$\tau$ is ergodic with respect to~$\mu$.
\end{lemma}

In turn, proving ergodicity reduces to verifying that a certain sequence of functionals converges to a constant in measure, and not to the more difficult problem of finding weak limits, as in the general ergodic method.

An illustrative example of the usefulness of the method of perturbations is given by the problem about central measures on the Young graph. Thoma's theorem about characters, more exactly, a rephrasing of this theorem as an assertion about the absolute of the Young graph, leaves no doubt that this problem should be regarded in connection with de Finetti's theorem. The fact is that the answers to these problems are remarkably similar. Namely, the absolute is stratified: a stratum of discrete measures parametrized by one-dimensional frequencies with sum~$1$ and a stratum of measures with zero frequencies. However, all proofs known to date (see~\cite{V3}) are not elementary and do not reveal the closeness of these problems. This relation is indeed nontrivial, and the main role here  is played by the dynamic properties of the RSK algorithm, which makes it possible to construct a required lifting of the graph of $Q$-tableaux to the Schur--Weyl graph.

Using $Q$-tableaux arising in the RSK algorithm to cover central measures on the Young graph by Bernoulli measures
was first suggested in~\cite{KV}; the fact that this correspondence is an isomorphism was proved in~\cite{SR}. But the embedding discussed above has not been noted; at the same time, a~careful analysis shows that the  method of perturbations allows one also to prove the theorem about the absolute, i.e., to prove that in this way we obtain all ergodic central measures. For Young tableaux with finitely many rows, this result is implicitly contained in~\cite{VT}.

\end{document}